\documentclass[12pt]{elsarticle}
\usepackage{mathrsfs}
\usepackage{amssymb,amsmath}
\usepackage{geometry}
\usepackage[colorlinks,
            linkcolor=red,
            anchorcolor=blue,
            citecolor=green
            ]{hyperref}
\newtheorem{theorem}{Theorem}[section]

\newtheorem{lemma}{Lemma}[section]

\newtheorem{remark}{Remark}[section]

\newcommand{\R}{\mathbb{R}}
\newcommand{\bal}{\begin{align}}
\newcommand{\bbal}{\begin{align*}}
\def\div{\mathord{{\rm div}}}

\newcommand{\De}{\Delta}
\newcommand{\na}{\nabla}

\newcommand{\cd}{\cdot}
\newcommand{\dd}{\mathrm{d}}

\newcommand{\pa}{\partial}

\newcommand{\fr}{\frac}
\newcommand{\ep}{\varepsilon}
\begin{document}
\begin{frontmatter}

\title{A class large solution of the 2D MHD equations with velocity and magnetic damping}

\author{Jinlu Li$^{a}$\quad Minghua Yang$^{b}$ and Yanghai Yu$^{c,}$\corref{cor1}}

\cortext[cor1]{Corresponding author. \\ \indent E-mail addresses: lijinlu@gnnu.cn (J. Li); ymh20062007@163.com (M. Yang); yuyanghai214@sina.com (Y. Yu);}
\address{$^{a}$ School of Mathematics and Computer Sciences, Gannan Normal University, Ganzhou 341000, China}
\address{$^{b}$ Department of Mathematics, Jiangxi University of Finance and Economics, Nanchang, 330032, China}
\address{$^{c}$ School of Mathematics and Statistics, Anhui Normal University, Wuhu, Anhui, 241002, China}

\begin{abstract}
  In this paper, we construct a class global large solution to the two-dimensional MHD equations with damp terms in the nonhomogeneous Sobolev framework.
  \end{abstract}


\end{frontmatter}

\section{Introduction}\label{sec1}
This paper focuses on the following 2D incompressible magnetohydrodynamics (MHD) equations
\begin{eqnarray}\label{2D-mhd}
        \left\{\begin{array}{ll}
          \partial_tu+u\cd\na u+\mu(-\Delta)^{\alpha} u+\na p=b\cd\na b,& x\in \R^2,t>0,\\
          \partial_tb+u\cd\na b+\nu(-\Delta)^{\beta} b=b\cd\na u,& x\in \R^2,t>0,\\
         \div u=\div b=0,& x\in \R^2,t\geq0,\\
          (u,b)|_{t=0}=(u_0,b_0),& x\in \R^2,\end{array}\right.
        \end{eqnarray}
where $u=(u_1(t,x),u_2(t,x))\in\R^2$ and $b=(b_1(t,x),b_2(t,x))\in\R^2$ denote the divergence free velocity field and magnetic field, respectively, $p\in \R$ is the scalar pressure. $\mu$ is the viscosity and $\nu$ is the magnetic diffusivity. The fractional power operator $(-\Delta)^{\gamma}$ with $0<\gamma<1$ is defined by Fourier multiplier with symbol $|\xi|^{2\gamma}$ (see e.g. \cite{Jacob 2005,Wu 2017})
\begin{eqnarray*}
   (-\Delta)^{\gamma}u(x)=\mathcal{F}^{-1}|\xi|^{2\gamma}\mathcal{F}u(\xi).
\end{eqnarray*}
We make the convention that by $\gamma= 0$ we mean that $(-\Delta)^{\gamma}u$ is a damp term $u$. The magnetohydrodynamic (MHD) equations which can be view as a coupling of incompressible Navier--Stokes and Maxwell's equations govern the motion of electrically conducting fluids such as plasmas, liquid metals and electrolytes, and play a fundamental role in geophysics, astrophysics, cosmology and engineering (see e.g. \cite{Priest 2000,Davidson 2001,Li 2017}). Due to the profound physical background and important mathematical significance, the MHD equations attracted quite a lot of
attention from many physicists and mathematicians in the past few years. Let us review some progress has been made about the MHD equations \eqref{2D-mhd} which are more relatively with our problem. It is well known that the 2D MHD equations \eqref{2D-mhd} with $-\Delta u$ and $-\Delta b$ (namely, $\alpha=\beta=1$) have the global smooth solution(\cite{Sermange}). In the completely inviscid case ($\mu=\nu=0$), the question of whether smooth solution of the MHD equations \eqref{2D-mhd} with large initial data develops singularity in finite time remains completely open. Besides these the two extreme cases, many intermediate cases, for example, the 2D MHD equations with partial dissipation, has been studied by various authors. The issue of the global regularity for the MHD equations \eqref{2D-mhd} with $\mu>0,\nu>0,\alpha>0,\beta=1$ has been solved by Fan et al.\cite{Fan 2014}. Recently, Yuan and Zhao \cite{Yuan 2018} considered the MHD equations \eqref{2D-mhd} with the dissipative operators weaker than any power of the fractional
Laplacian and obtained the global regularity of the corresponding system. On the other hand, Cao et al.\cite{Cao 2014}, Jiu and Zhao \cite{Jiu 2015} established the global regularity of smooth solutions to the MHD equations \eqref{2D-mhd} with $\mu=0,\nu>0,\beta>1$ by different approach. Subsequently, Agelas \cite{Agelas 2016} improved this work with the diffusion $(-\Delta)^{\beta} b (\beta > 1)$ replaced by $(-\Delta)\log^{\kappa}(e-\Delta) b (\kappa > 1)$.

As mentioned above, the global regularity for the completely inviscid MHD equations \eqref{2D-mhd} with large initial data is still a challenging open problem. When $\alpha=\beta=0$, Wu et al \cite{Wu 2015} obtained that the d-dimensional MHD equations \eqref{2D-mhd} always possesses a unique global solution provided that the initial datum is sufficiently small in the nonhomogeneous functional setting $H^s$ with $s>1+\fr d2$. Our main goal is to prove the global existence of solutions to \eqref{2D-mhd} with $\alpha=\beta=0$ for a class of large initial data.

We assume from now on that the damping coefficients $\mu=\nu=1$, just for simplicity. Our main result is stated as follows.
\begin{theorem}\label{the1.1} Let $\alpha=\beta=0$ and $s>2$. Assume that the initial data fulfills ${\rm{div}}u_0={\rm{div}}b_0=0$ and
$$u_0=U_0+v_0\quad \mbox{and}\quad b_0=B_0+c_0$$
where
\begin{eqnarray*}
&U_0=
\begin{pmatrix}
\pa_2a_0 \\ -\pa_1a_0
\end{pmatrix}
\quad\mbox{and}\quad
B_0=
\begin{pmatrix}
\pa_2m_0 \\ -\pa_1m_0
\end{pmatrix}
\end{eqnarray*}
with
\begin{eqnarray}\label{Equ1.2}
\mathrm{supp} \ \hat{a}_0(\xi),\mathrm{supp}\ \hat{m}_0(\xi)\subset\mathcal{C}:=\Big\{\xi \big| \ |\xi_1-\xi_2|\leq \ep\Big\} .
\end{eqnarray}
There exists a sufficiently small positive constant $\delta$, and a universal constant $C$ such that if
\begin{align}\label{condition}
\Big(||v_0||^2_{H^s}+||c_0||^2_{H^s}+\ep^2(||a_0||^4_{H^s}+||m_0||^4_{H^s})\Big)
\exp\Big(C(||a_0||_{H^{s+2}}+||m_0||_{H^{s+2}})\Big)\leq \delta,
\end{align}
then the system \eqref{2D-mhd} has a unique global solution.
\end{theorem}
\begin{remark}\label{rem1.1}
Let $v_0=c_0=0$ and $a_0=m_0=\ep^{-\frac12}\log\log\frac1\ep \chi$, where the smooth function $\chi$ satisfying
\begin{align*}
\mathrm{supp} \hat{\chi}\in \mathcal{\widetilde{C}},\quad \hat{\chi}(\xi)\in[0,1]\quad\mbox{and} \quad \hat{\chi}(\xi)=1 \quad\mbox{for} \quad \xi\in\frac12\mathcal{\widetilde{C}},
\end{align*}
where
\begin{align*}
\mathcal{\widetilde{C}}:=\Big\{\xi \big| \ |\xi_1-\xi_2|\leq \ep,\ 1\leq\xi^2_1+\xi^2_2\leq 2\Big\} .
\end{align*}
Then, direct calculations show that the left side of \eqref{condition} becomes
\begin{align*}
C\ep^2\Big(\log\log \frac1\ep\Big)^4\exp\Big(C\log\log \frac1\ep\Big).
\end{align*}
Therefore, choosing $\ep$ small enough, we deduce that the system \eqref{2D-mhd} has a global solution.

Moreover, we also have
\begin{align*}
||u_0||_{L^2}\gtrsim \log\log \frac1\ep\quad\mbox{and}\quad||b_0||_{L^2}\gtrsim \log\log \frac1\ep.
\end{align*}
\end{remark}
\begin{remark}\label{rem1.2} Considered the system \eqref{2D-mhd} with $0<\alpha,\beta<1$, if the support condition \eqref{Equ1.2} of the Theorem \ref{the1.1} were replaced by
\begin{eqnarray}\label{l}
\mathrm{supp} \ \hat{a}_0(\xi),\mathrm{supp}\ \hat{m}_0(\xi)\subset\mathcal{C}:=\Big\{\xi \big| \ |\xi_1-\xi_2|\leq \ep,\ 1\leq\xi^2_1+\xi^2_2\Big\},
\end{eqnarray}
the Theorem \ref{the1.1} holds true.
\end{remark}
{\bf Notations}: For the sake of simplicity, $a\lesssim b$ means that there is a uniform positive constant $C$ such that $a\leq Cb$. $[A,B]$ stands for the commutator operator $AB-BA$, where $A$ and $B$ are any pair of operators on some Banach space. In the paper, we will use the Besov space ${B}_{p,q}^{s}$, for more details, we refer the readers to see the Chapter 2 in \cite{Bahouri2011}. It is worth mentioning that the Besov space ${B}_{2,2}^{s}$ coincides with the nonhomogeneous Sobolev spaces $H^s$ for $s>0$, namely, ${{B}_{2,2}^{s}}(\R^{d})={H^s}(\R^{d}),$
where
 $$H^s(\R^{d}):=\Big\{\mathbf{f}\in \mathcal{S'}(\R^{d}):||\mathbf{f}||_{H^s(\R^{d})}<\infty\Big\}$$
 with the norm
$$
||\mathbf{f}||_{H^s(\R^{d})}:=\Big(\int_{\R^{d}}(1+|\xi|^2)^s|\widehat{f}(\xi)|^2d\xi\Big)^{\fr12}.$$
\section{Reformulation of the System}\label{sec2}
\setcounter{equation}{0}
Let $(a,m)$ be the solutions of the following system
\begin{eqnarray}\label{l1}
        \left\{\begin{array}{ll}
          \pa_ta+a=0,\\
          \pa_tm+m=0,\\
          (a,m)|_{t=0}=(a_0,m_{0}).\end{array}\right.
        \end{eqnarray}
Setting
\bbal
&U=
\begin{pmatrix}
\pa_2a \\ -\pa_1a
\end{pmatrix}
 \quad\mbox{and}\quad
B=
\begin{pmatrix}
\pa_2m \\ -\pa_1m
\end{pmatrix},
\end{align*}
we can deduce from \eqref{l1} that
\begin{eqnarray}\label{l2}
        \left\{\begin{array}{ll}
          \pa_t U+U=0,\\
          \pa_t B+B=0,\\
          \div U=\div B=0,\\
          (U,B)|_{t=0}=(U_0,B_{0}).\end{array}\right.
        \end{eqnarray}
Denoting $v=u-U$ and $c=b-B$, the system \eqref{2D-mhd} can be written as follows
\begin{eqnarray}\label{l3}
        \left\{\begin{array}{ll}
\partial_tv+u\cd\na v+v\cd\na U+u+\na p=b\cd\na c+c\cd\na B+f,\\
\partial_tc+u\cd\na c+v\cd\na B+c=b\cd\na v+c\cd\na U+g,\\
\div v=\div c=0,\\
(v,c)|_{t=0}=(v_0,c_0).\end{array}\right.
\end{eqnarray}
where
\bbal
&f=-U\cd\na U+B\cd\na B\quad \mbox{and}\quad g=-U\cd\na B+B\cd\na U.
\end{align*}
\section{The Proof of Theorem \ref{the1.1}}\label{sec3}
\setcounter{equation}{0}
Before proceeding on, we present some estimates which will be used in the proof of Theorem \ref{the1.1}.

\begin{lemma}\label{lem3.1} For $s>2$, under the assumptions of Theorem \ref{the1.1}, the following estimates hold
\bal\label{l}
||f||_{H^s}+||g||_{H^s}\leq Ce^{-t}\ep(||a_0||^2_{H^{s+2}}+||m_0||^2_{H^{s+2}})
\end{align}
and
\bal\label{lj}
||\na U||_{H^s}+||\na B||_{H^s}\leq Ce^{-t}(||a_0||_{H^{s+2}}+||m_0||_{H^{s+2}}).
\end{align}
\end{lemma}
{\bf Proof of Lemma \ref{lem3.1}}\quad Notice that
\bbal
f^1&=-U\cd\na U^1+B\cd\na B^1
\\&=\pa_1a\pa_2\pa_2a-\pa_2a\pa_1\pa_2a-(\pa_1m\pa_2\pa_2m-\pa_2m\pa_1\pa_2m)
\\&=(\pa_1-\pa_2)a\pa_2\pa_2a+\pa_2a\pa_2(\pa_2-\pa_1)a+(\pa_2-\pa_1)m\pa_2\pa_2m+\pa_2m\pa_2(\pa_1-\pa_2)m
\end{align*}
and
\bbal
f^2&=-U\cd\na U^2+B\cd\na B^2
\\&=-\pa_1a\pa_2\pa_1a+\pa_2a\pa_1\pa_1a+\pa_1m\pa_2\pa_1m-\pa_2m\pa_1\pa_1m
\\&=(\pa_2-\pa_1)a\pa_1\pa_2a+\pa_2a\pa_1(\pa_1-\pa_2)a+(\pa_1-\pa_2)m\pa_1\pa_2m+\pa_2m\pa_1(\pa_2-\pa_1)m,
\end{align*}
due to the fact that $H^s$ with $s>2$ is a Banach algebra, then we have
\bal\label{lj1}
||f^1||_{H^s}\lesssim&~ ||(\pa_1-\pa_2)a||_{H^s}||a||_{H^{s+2}}+||a||_{H^{s+1}}||(\pa_2-\pa_1)a||_{H^{s+1}}\nonumber\\
&+||(\pa_1-\pa_2)m||_{H^s}||m||_{H^{s+2}}+||m||_{H^{s+1}}||(\pa_2-\pa_1)m||_{H^{s+1}}.
\end{align}
Direct calculations show that for $\tau\geq0$
\bal\label{lj2}
||a||_{H^{\tau}}+||m||_{H^{\tau}}\leq e^{-t}(||a_0||_{H^{\tau}}+||m_0||_{H^{\tau}})
\end{align}
and
\bal\label{lj3}
&||(\pa_1-\pa_2)a||_{H^{\tau}}+||(\pa_1-\pa_2)m||_{H^{\tau}}\nonumber\\
\leq&~ e^{-t}(||(\pa_1-\pa_2)a_0||_{H^{\tau}}+||(\pa_1-\pa_2)a_0||_{H^{\tau}})\nonumber\\
\leq&~ e^{-t}\ep(||a_0||_{H^{\tau}}+||m_0||_{H^{\tau}}),
\end{align}
where we have used the conditions $\mathrm{supp}\ \hat{a}_0(\xi)\subset \mathcal{C}$ and $\mathrm{supp}\ \hat{m}_0(\xi)\subset \mathcal{C}$.

In view of the facts \eqref{lj2} and \eqref{lj3}, we obtain from \eqref{lj1} that
\bbal
||f^1||_{H^s}\lesssim&~e^{-t}\ep(||a_0||^2_{H^{s+2}}+||m_0||^2_{H^{s+2}}).
\end{align*}
Similarly, we also have
\bbal
||f^2||_{H^s}\lesssim&~e^{-t}\ep(||a_0||^2_{H^{s+2}}+||m_0||^2_{H^{s+2}}).
\end{align*}
Then, we get
\bal\label{lj4}
||f||_{H^s}\leq||f^1||_{H^s}+||f^2||_{H^s}\lesssim&~e^{-t}\ep(||a_0||^2_{H^{s+2}}+||m_0||^2_{H^{s+2}}).
\end{align}
An argument similar to that used above, we get
\bal\label{lj5}
||g||_{H^s}\lesssim&~e^{-t}\ep(||a_0||^2_{H^{s+2}}+||m_0||^2_{H^{s+2}}).
\end{align}
Combining \eqref{lj4} and \eqref{lj5} yields the desired result \eqref{l}.

\eqref{lj} is just a consequence of \eqref{lj2}. Thus, we complete the proof of Lemma \ref{lem3.1}. $\Box$

{\bf Proof of Theorem \ref{the1.1}}\quad For notational simplicity, we set
\bbal
&E(t)=\big(||v(t)||^2_{H^s}+||c(t)||^2_{H^s}\big).
\end{align*}
Applying $\De_j$ to \eqref{l3} and taking the $L^2$ inner product of the resulting equations with $\De_jv$ and $\De_jc$, respectively, we have
\bal\label{j}
&\frac12\frac{\dd}{\dd t}\Big(||\De_jv||^2_{L^2}+||\De_jc||^2_{L^2}\Big)+||\De_jv||^2_{L^2}
+||\De_jc||^2_{L^2}=:\sum_{i=1}^5K_i,
\end{align}
where
\bbal
&K_1=-\int_{\R^2}[\De_j,u\cd\na]v\cd\De_jv\dd x-\int_{\R^2}[\De_j,u\cd\na]c\cd\De_jc\dd x,
\\&K_2=\int_{\R^2}[\De_j,b\cd\na]c\cd\De_jv\dd x+\int_{\R^2}[\De_j,b\cd\na]v\cd\De_jc\dd x,
\\&K_3=-\int_{\R^2}\De_j(v\cd\na U)\cd\De_jv\dd x-\int_{\R^2}\De_j(v\cd\na  B)\cd\De_jc\dd x,
\\&K_4=\int_{\R^2}\De_j(c\cd\na B)\cd\De_jv\dd x+\int_{\R^2}\De_j(c\cd\na U)\cd\De_jc\dd x,
\\&K_5=\int_{\R^2}\De_jf\cd\De_jv\dd x+\int_{\R^2}\De_jg\cd\De_jc\dd x.
\end{align*}

Multiplying both sides of \eqref{j} by $2^{2js}$ and summing up over $j\geq-1$ yields
\bal\label{j1}
&\frac12\frac{\dd}{\dd t}E(t)+E(t)=\sum_{i=1}^5\sum_{j\geq -1}2^{2js}K_i.
\end{align}

Next, we need to estimate the above terms involving $K_i$ for $i=1,\cdot\cdot\cdot,5$ as follows
\bal\label{y1}
\sum_{j\geq -1}2^{2js}|K_1|\leq&~\sum_{j\geq -1}2^{2js}||[\De_j,u\cd\na]v||_{L^2}||\De_jv||_{L^2}
+\sum_{j\geq -1}2^{2js}||[\De_j,u\cd\na]c||_{L^2}||\De_jc||_{L^2} \nonumber\\
\lesssim&~||\na u||_{H^{s-1}}\Big(||v||^2_{H^s}+||c||^2_{H^s}\Big)\nonumber\\
\lesssim&~\Big(||U||_{H^s}+||v||_{H^s}\Big)E(t),
\end{align}
where we have used the commutator estimate (see Lemma 2.6 in \cite{Wu 2015})
$$||[\De_j,u\cd\na]\mathbf{f}||_{B^{s}_{2,2}}\leq C||\na u||_{B^{s-1}_{2,2}}||\mathbf{f}||_{B^{s}_{2,2}}\quad\mbox{with}\quad\div u=0.$$

Similarly, we also have
\bal\label{y2}
\sum_{j\geq -1}2^{2js}|K_2|&\lesssim\Big(||B||_{H^s}+||c||_{H^s}\Big)E(t).
\end{align}

For the last three terms, by H\"{o}lder's inequality, we deduce
\bal
\sum_{j\geq -1}2^{2js}|K_3|&\lesssim ||\na U||_{H^s}||v||^2_{H^s}+||\na B||_{H^s}||c||_{H^s}||v||_{H^s},\label{y3}\\
\sum_{j\geq -1}2^{2js}|K_4|&\lesssim ||\na U||_{H^s}||c||^2_{H^s}+||\na B||_{H^s}||c||_{H^s}||v||_{H^s},\label{y4}
\end{align}
and
\bal\label{y5}
\sum_{j\geq -1}2^{2js}|K_5|&\leq \sum_{j\geq -1}2^{2js}||\De_jf||_{L^2}||\De_jv||_{L^2}+\sum_{j\geq -1}2^{2js}||\De_jg||_{L^2}||\De_jc||_{L^2}\nonumber\\
&\leq C(||f||^2_{H^s}+||g||^2_{H^s})+\frac12E(t).
\end{align}
Inserting \eqref{y1}--\eqref{y5} into \eqref{j1} yields that
\bal\label{yy4}
\frac{\dd}{\dd t}E(t)+E(t)&\lesssim E^{\frac32}(t)+\Big(||\na U||_{H^s}+||\na B||_{H^s}\Big)E(t)+||f||^2_{H^s}+||g||^2_{H^s}.
\end{align}
Utilizing the Lemma \ref{lem3.1}, we have from \eqref{yy4}
\bal\label{yy5}
\frac{\dd}{\dd t}E(t)+E(t)\lesssim&~ E^{\frac32}(t)+e^{-t}\Big(||a_0||_{H^{s+2}}+||m_0||_{H^{s+2}}\Big)E(t)+e^{-t}\ep^2\Big(||a_0||^4_{H^{s+2}}+||m_0||^4_{H^{s+2}}\Big).
\end{align}

Now, we define
\bbal
\Gamma:=\max\{t\in[0,T^*):\sup_{\tau\in[0,t]}E(\tau)\leq \eta\},
\end{align*}
where $\eta$ is a small enough positive constant which will be determined later on.

Assume that $\Gamma<T^*$. For all $t\in[0,\Gamma]$, we obtain from \eqref{yy5} that
\bbal
\frac{\dd}{\dd t}E(t) \leq Ce^{-t}\Big(||a_0||_{H^{s+2}}+||m_0||_{H^{s+2}}\Big)E(t)
+Ce^{-t}\ep^2\Big(||a_0||^4_{H^{s+2}}+||m_0||^4_{H^{s+2}}\Big),
\end{align*}
which follows from the assumption \eqref{condition} that
\bbal
E(t)&\leq C\Big(E_0+\ep^2(||a_0||^4_{H^s}+||m_0||^4_{H^s})\Big)\exp\Big(C(||a_0||_{H^{s+2}}+||m_0||_{H^{s+2}})\Big)\leq C\delta.
\end{align*}
Choosing $\eta=2C\delta$, thus we can get
\bbal
E(t)&\leq \fr\eta2 \quad\mbox{for}\quad t\leq \Gamma.
\end{align*}

So if $\Gamma<T^*$, due to the continuity of the solutions, we can obtain there exists $0<\epsilon\ll1$ such that
\bbal
E(t)&\leq \fr\eta2 \quad\mbox{for}\quad t\leq \Gamma+\epsilon<T^*,
\end{align*}
which is contradiction with the definition of $\Gamma$.

Thus, we can conclude $\Gamma=T^*$ and
\bbal
E(t)&\leq C<\infty \quad\mbox{for all}\quad t\in(0,T^*),
\end{align*}
which implies that $T^*=+\infty$. This completes the proof of Theorem \ref{the1.1}. $\Box$

\section*{Acknowledgments} J. Li was partially supported by NSFC (No.11801090). M. Yang was partially supported by NSFC (No.11801236)


\section*{References}
\begin{thebibliography}{00}\label{ref:ref}\addtolength{\itemsep}{-1.6ex}
\bibitem{Agelas 2016} L. Agelas, Global regularity for logarithmically critical 2D MHD equations with zero viscosity. Monatsh. Math. 181 (2016), 245--266.
\bibitem{Bahouri2011} H. Bahouri, J.Y. Chemin, R. Danchin, Fourier Analysis and Nonlinear Partial Differential Equations, Grundlehren Math. Wiss., vol.343, Springer-Verlag, Berlin, Heidelberg, 2011.
\bibitem{Cao 2014} C. Cao, J. Wu, B. Yuan, The 2D incompressible magnetohydrodynamics equations with only magnetic diffusion. SIAM J. Math. Anal. 46 (2014), 588--602.
\bibitem{Davidson 2001} P.A. Davidson, An Introduction to Magnetohydrodynamics, Cambridge University Press, Cambridge, England, 2001.
\bibitem{Fan 2014} J. Fan, H. Malaikah, S. Monaquel, G. Nakamura, Y. Zhou, Global cauchy problem of 2D generalized MHD equations. Monatsh Math. 175  (2014), 127--131.
\bibitem{Sermange} M. Sermange, R. Temam, Some mathematical questions related to the MHD equations. Comm. Pure Appl. Math. 36 (1983), 635--664.
\bibitem{Jacob 2005} N. Jacob, Pseudo Differential Operators and Markov Processes. Vol. III: Markov Processes and Applications, Imperial College Press, 2005.
\bibitem{Jiu 2015} Q. Jiu, J. Zhao, Global regularity of 2D generalized MHD equations with magnetic diffusion. Z. Angew. Math. Phys. 66 (2015), 677--687.
\bibitem{Li 2017} J. Li, W. Tan, Z. Yin, Local existence and uniqueness for the non-resistive MHD equations in homogeneous Besov spaces. Advances in Mathematics. 317 (2017) 786--798.
\bibitem{Priest 2000} E. Priest, T. Forbes, Magnetic Reconnection, MHD Theory and Applications, Cambridge University Press, Cambridge, 2000.
\bibitem{Wu 2015} J. Wu , X. Xu, Z. Ye, Global smooth solutions to the n-dimensional damped models of incompressible fluid mechanics with small initial datum. J. Nonlinear Sci. 25 (2015), 157--192.
\bibitem{Wu 2017} X. Wu, Y. Yu, Y. Tang, Global existence and asymptotic behavior for the 3D generalized Hall-MHD system. Nonlinear Anal. 151 (2017) 41--50.
\bibitem{Yuan 2018} B. Yuan, J. Zhao, Global regularity of 2D almost resistive MHD equations. Nonlinear Anal. Real World Appl. 41 (2018), 53--65.


\end{thebibliography}
\end{document}